\documentclass[12pt,twoside]{article}
\usepackage[english]{babel}
\usepackage[latin1]{inputenc}
\usepackage{amsmath}
\usepackage{amssymb,amsfonts}
\usepackage{graphicx}
\usepackage{times,amssymb,amscd}

\newcommand{\bC}{\mathbf{C}}
\newcommand{\bE}{\mathbf{E}}
\newcommand{\bG}{\mathbf{G}}
\newcommand{\bH}{\mathbf{H}}

\newcommand{\bL}{\mathbf{L}}

\newcommand{\bR}{\mathbf{R}}
\newcommand{\bS}{\mathbf{S}}

\newcommand{\bs}{\mathbf{s}}
\newcommand{\ba}{\mathbf{a}}

\newcommand{\bT}{\mathbf{T}}

\newcommand{\bb}{\mathbf{b}}

\newcommand{\cP}{\mathcal{P}}
\newcommand{\cS}{\mathcal{S}}
\newcommand{\cT}{\mathcal{T}}

\newcommand{\cH}{\mathcal{H}}

\newcommand{\EUC}{\bE^3}

\newcommand{\SXR}{\bS^2\!\times\!\bR}
\newcommand{\HXR}{\bH^2\!\times\!\bR}
\newcommand{\SLR}{\widetilde{\bS\bL_2\bR}}
\newcommand{\NIL}{\mathbf{Nil}}
\newcommand{\SOL}{\mathbf{Sol}}

\newtheorem{Definition}{Definition}[section]
\newtheorem{Remark}{Remark}[section]

\begin{document}
\pagestyle{myheadings}
\markboth{\centerline{Jen\H o Szirmai}}
{Non-periodic geodesic ball packings...}
\title
{Non-periodic geodesic ball packings to infinite regular prism tilings in $\SLR$ space 
}

\author{Jen\H o Szirmai \\
\normalsize Budapest University of Technology and \\
\normalsize Economics Institute of Mathematics, \\
\normalsize Department of Geometry \\
\normalsize Budapest, P. O. Box: 91, H-1521 \\
\normalsize szirmai@math.bme.hu
\date{\normalsize{\today}}}

%%%%%%%%%%%%%%%%%%%%%%%%%%%%%%%%%%%%%%%%%%%%
%%AMS Classification 2000
%% The 1-st classification is obligatory, the 2-nd classification is optional
%% \subjclass{primary}{secondary}      f.e. \subjclass{35R35, 49N50}{}
%%{52C17, 52C22}{}
%%%%%%%%%%%%%%%%%%%%%%%%%%%%%%%%%%%%%%%%%%%%

\maketitle
\begin{abstract}

In \cite{Sz13-1} we defined and described the {\it regular infinite or bounded} $p$-gonal prism tilings in $\SLR$ space.
We proved that there exist infinitely many regular infinite $p$-gonal face-to-face prism tilings $\cT^i_p(q)$ and
infinitely many regular bounded $p$-gonal non-face-to-face prism tilings $\cT_p(q)$ for integer parameters $p,q;~3 \le p$,
$ \frac{2p}{p-2} < q$. Moreover, in \cite{MSz14} and \cite{MSzV13} we have determined the symmetry group of $\cT_p(q)$ 
via its index 2 rotational subgroup, denoted by $\mathbf{pq2_1}$ and investigated 
the corresponding geodesic and translation ball packings. 

In this paper we study the structure of the regular infinite or bounded $p$-gonal prism tilings, prove that the side curves of their base figurs
are arcs of Euclidean circles  for each parameter. Moreover, we examine the non-periodic geodesic  
ball packings of congruent regular non-periodic prism tilings derived from the regular infinite $p$-gonal 
face-to-face prism tilings $\cT^i_p(q)$ in $\SLR$ geometry. We develop a procedure to determine
the densities of the above non-periodic optimal geodesic ball packings and apply this algorithm to them.
We look for those parameters $p$ and $q$ above, where the packing density large enough as possible. 
Now, we obtain larger density $\approx 0.626606$ for $(p, q) = (29,3)$ 
then the maximal density of the corresponding periodical geodesic ball packings under the groups $\mathbf{pq2_1}$. 

In our work we will use the projective model of $\SLR$ introduced by E. {Moln\'ar} in \cite{M97}.
\end{abstract}
%%%%%%%%%%%%%%%%%%%%%%%%%%%%%%%%%%%%%%%%%%%%

%%%%%%%%%%%%%%%%%%%%%%%%%%%%%%%%%%%%%%%%%%%
\newtheorem{theorem}{Theorem}[section]
\newtheorem{corollary}[theorem]{Corollary}
\newtheorem{conjecture}[theorem]{Conjecture}
\newtheorem{lemma}[theorem]{Lemma}
\newtheorem{exmple}[theorem]{Example}
\newenvironment{definition}{\begin{defn}\normalfont}{\end{defn}}
\newenvironment{example}{\begin{exmple}\normalfont}{\end{exmple}}
\newenvironment{acknowledgement}{Acknowledgement}

%%%%%%%%%%%%%%%%%%%%%%%%%%%%%%%%%%%%%%%%%%%%%%%%%%%%%%%%%%%%%%%%%%%%

%============================================================================%
%                             the main article                               %
%============================================================================%

%%%%%%%%%%%%%%%%%%%%%%%%%%%%%%%%%%%%%%%%%%%%%%%%%%%%%%%%%%%%%%%%%%%%%%%%%%%%%%
\section{Basic notions}

The real $ 2\times 2$ matrices $\begin{pmatrix}
         d&b \\
         c&a \\
         \end{pmatrix}$ with unit determinant $ad-bc=1$
constitute a Lie transformation group by the usual product operation, taken to act on row matrices as on point coordinates on the right as follows
\begin{equation}
\begin{gathered}
(z^0,z^1)\begin{pmatrix}
         d&b \\
         c&a \\
         \end{pmatrix}=(z^0d+z^1c, z^0 b+z^1a)=(w^0,w^1)\\
\mathrm{with} \ w=\frac{w^1}{w^0}=\frac{b+\frac{z^1}{z^0}a}{d+\frac{z^1}{z^0}c}=\frac{b+za}{d+zc}, \tag{1.1}
\end{gathered}
\end{equation}
as action on the complex projective line $\bC^{\infty}$ (see \cite{M97}, \cite{MSz}).
This group is a $3$-dimensional manifold, because of its $3$ independent real coordinates and with its usual neighbourhood topology (\cite{S}, \cite{T}, \cite{R}).
In order to model the above structure in the projective sphere $\cP \cS^3$ and in the projective space $\cP^3$ (see \cite{M97}),
we introduce the new projective coordinates $(x^0,x^1,x^2,x^3)$ where
$a:=x^0+x^3, \ b:=x^1+x^2, \ c:=-x^1+x^2, \ d:=x^0-x^3$
with the positive, then the non-zero multiplicative equivalence as projective freedom in $\cP \cS^3$ and in $\cP^3$, respectively.
Then it follows that $0>bc-ad=-x^0x^0-x^1x^1+x^2x^2+x^3x^3$
describes the interior of the above one-sheeted hyperboloid solid $\cH$ in the usual Euclidean coordinate simplex with the origin
$E_0(1;0;0;0)$ and the ideal points of the axes $E_1^\infty(0;1;0;0)$, $E_2^\infty(0;0;1;0)$, $E_3^\infty(0;0;0;1)$.
We consider the collineation group ${\bf G}_*$ that acts on the projective sphere $\cS\cP^3$  and preserves a polarity i.e. a scalar product of signature
$(- - + +)$, this group leaves the one sheeted hyperboloid solid $\cH$ invariant.
We have to choose an appropriate subgroup $\mathbf{G}$ of $\mathbf{G}_*$ as isometry group, then the universal covering group and space
$\widetilde{\cH}$ of $\cH$ will be the hyperboloid model of $\SLR$ \cite{M97}.

The specific isometries $\bS(\phi)$ $(\phi \in \bR )$ constitute a one parameter group given by the matrices:
\begin{equation}
\begin{gathered} \bS(\phi):~(s_i^j(\phi))=
\begin{pmatrix}
\cos{\phi}&\sin{\phi}&0&0 \\
-\sin{\phi}&\cos{\phi}&0&0 \\
0&0&\cos{\phi}&-\sin{\phi} \\
0&0&\sin{\phi}&\cos{\phi}
\end{pmatrix}
\end{gathered} \tag{1.2}
\end{equation}
The elements of $\bS(\phi)$ are the so-called {\it fibre translations}. We obtain a unique fibre line to each $X(x^0;x^1;x^2;x^3) \in \widetilde{\cH}$
as the orbit by right action of $\bS(\phi)$ on $X$. The coordinates of points lying on the fibre line through $X$ can be expressed
as the images of $X$ by $\bS(\phi)$:
\begin{equation}
\begin{gathered}
(x^0;x^1;x^2;x^3) \stackrel{\bS(\phi)}{\longrightarrow} {(x^0 \cos{\phi}-x^1 \sin{\phi}; x^0 \sin{\phi} + x^1 \cos{\phi};} \\ {x^2 \cos{\phi} + x^3 \sin{\phi};-x^2 \sin{\phi}+
x^3 \cos{\phi})}.
\end{gathered} \tag{1.3}
\end{equation}
The points of a fibre line through
$X$ by usual inhomogeneous Euclidean coordinates $x=\frac{x^1}{x^0}$, $y=\frac{x^2}{x^0}$, $z=\frac{x^3}{x^0}$, $x^0\ne 0$ are given by
\begin{equation}
\begin{gathered}
(1;x;y;z) \stackrel{\bS(\phi)}{\longrightarrow} {\Big( 1; \frac{x+\tan{\phi}}{1-x \tan{\phi}}; \frac{y+z \tan{\phi}}{1-x \tan{\phi}};
\frac{z - y \tan{\phi}}{1-x \tan{\phi}}\Big)}
\end{gathered} \tag{1.4}
\end{equation}
for the projective space $\cP^3$, where ideal points (at infinity) conventionally occur.

In (1.3) and (1.4) we can see the $2\pi$ periodicity of $\phi$, moreover the (logical) extension to $\phi \in \bR$, as real parameter, to have
the universal covers $\widetilde{\cH}$ and $\SLR$, respectively, through the projective sphere $\cP\cS^3$. The elements of the isometry group of
$\mathbf{SL_2R}$ (and so by the above extension the isometries of $\SLR$) can be described by the matrix $(a_i^j)$ (see \cite{M97} and \cite{MSz})
Moreover, we have the projective proportionality, of course.
We define the {\it translation group} $\bG_T$, as a subgroup of the isometry group of $\mathbf{SL_2R}$,
the isometries acting transitively on the points of ${\cH}$ and by the above extension on the points of $\SLR$ and $\widetilde{\cH}$.
$\bG_T$ maps the origin $E_0(1;0;0;0)$ onto $X(x^0;x^1;x^2;x^3)$. These isometries and their inverses (up to a positive determinant factor) 
can be given by the following matrices:
\begin{equation}
\begin{gathered} \bT:~(t_i^j)=
\begin{pmatrix}
x^0&x^1&x^2&x^3 \\
-x^1&x^0&x^3&-x^2 \\
x^2&x^3&x^0&x^1 \\
x^3&-x^2&-x^1&x^0
\end{pmatrix}.
\end{gathered} \tag{1.5}
\end{equation}
The rotation about the fibre line through the origin $E_0(1;0;0;0)$ by angle $\omega$ $(-\pi<\omega\le \pi)$ can be expressed by the following matrix
(see \cite{M97})
\begin{equation}
\begin{gathered} \bR_{E_O}(\omega):~(r_i^j(E_0,\omega))=
\begin{pmatrix}
1&0&0&0 \\
0&1&0&0 \\
0&0&\cos{\omega}&\sin{\omega} \\
0&0&-\sin{\omega}&\cos{\omega}
\end{pmatrix},
\end{gathered} \tag{1.6}
\end{equation}
and the rotation $\bR_X(\omega)$ about the fibre line through $X(x^0;x^1;x^2;x^3)$ by angle $\omega$ can be derived by formulas (1.5) and (1.6):
\begin{equation}
\bR_X(\omega)=\bT^{-1} \bR_{E_O} (\omega) \bT:~(r_i^j(X,\omega)).
\tag{1.7}
\end{equation}
Horizontal intersection of the hyperboloid solid $\cH$ with the plane $E_0 E_2^\infty E_3^\infty$ provides the
{\it hyperbolic $\mathbf{H}^2$ base plane} of the model $\widetilde{\cH}=\SLR$.
The fibre through $X$ intersects the base plane $z^1=x=0$ in the foot point
\begin{equation}
\begin{gathered}
Z(z^0=x^0 x^0+x^1x^1; z^1=0; z^2=x^0x^2-x^1x^3;z^3=x^0x^3+x^1x^2).
\end{gathered} \tag{1.8}
\end{equation}
We introduce a so-called hyperboloid parametrization by \cite{M97} as follows
\begin{equation}
\begin{gathered}
x^0=\cosh{r} \cos{\phi}, ~ ~
x^1=\cosh{r} \sin{\phi}, \\
x^2=\sinh{r} \cos{(\theta-\phi)}, ~ ~
x^3=\sinh{r} \sin{(\theta-\phi)},
\end{gathered} \tag{1.9}
\end{equation}
where $(r,\theta)$ are the polar coordinates of the base plane and $\phi$ is just the fibre coordinate. We note that
$$-x^0x^0-x^1x^1+x^2x^2+x^3x^3=-\cosh^2{r}+\sinh^2{r}=-1<0.$$
The inhomogeneous coordinates corresponding to (1.9), that play an important role in the later visualization of prism tilings in $\EUC$,
are given by
\begin{equation}
\begin{gathered}
x=\frac{x^1}{x^0}=\tan{\phi}, ~ ~
y=\frac{x^2}{x^0}=\tanh{r} \frac{\cos{(\theta-\phi)}}{\cos{\phi}}, \\
z=\frac{x^3}{x^0}=\tanh{r} \frac{\sin{(\theta-\phi)}}{\cos{\phi}}.
\end{gathered} \tag{1.10}
\end{equation}
\subsection{Geodesic balls in $\SLR$}
\begin{Definition}
The {\rm distance} $d(P_1,P_2)$ between the points $P_1$ and $P_2$ is defined by the arc length of the geodesic curve
from $P_1$ to $P_2$.
\end{Definition}
\begin{Definition}
The {\rm geodesic sphere} of radius $\rho$ (denoted by $S_{P_1}(\rho)$) with the center in point $P_1$ is defined as the set of all points
$P_2$ with the condition $d(P_1,P_2)=\rho$. Moreover, we require that the geodesic sphere is a simply connected
surface without selfintersection.
\end{Definition}
\begin{Definition}
The body of the geodesic sphere of centre $P_1$ and with radius $\rho$ is called {\rm geodesic ball}, denoted by $B_{P_1}(\rho)$,
i.e., $Q \in B_{P_1}(\rho)$ iff $0 \leq d(P_1,Q) \leq \rho$.
\end{Definition}
From \cite{MSz14} it follows that $S(\rho)$ is a simply connected surface in $\mathbf{E}^3$ and $\SLR$, respectively, if $\rho \in [0,\frac{\pi}{2})$.
If $\rho\ge \frac{\pi}{2}$ then the universal cover should be discussed.
{\it Therefore, we consider geodesic spheres and balls only with radii $\rho \in [0,\frac{\pi}{2})$ in the following}.
\subsection{The volume of a geodesic ball}
The volume formula of the geodesic ball $B(\rho)$ follows from the metric tensor $g_{ij}$ (see \cite{MSz14}).
We obtain the connection between the hyperboloid coordinates $(r,\theta,\phi)$ and the geographical coordinates $(s,\lambda,\alpha)$
in a standard way. Therefore, the volume of the geodesic ball of radius $\rho$ can be computed by the
following
\begin{theorem}
\begin{equation}
\begin{gathered}
Vol(B(\rho))=\int_{B} \frac{1}{2}\sinh(2r)~ dr ~d \theta~ d \phi  = \\ = 4\pi \int_{0}^{\rho} \int_{0}^{\frac{\pi}{4}}
\frac{1}{2}\sinh(2r(s,  \alpha ))\dot |J_1|  ~ d \alpha \ ds
\\ +4 \pi \int_{0}^{\rho} \int_{\frac{\pi}{4}}^{\frac{\pi}{2}}
\frac{1}{2}\sinh(2r(s,  \alpha ))\dot |J_2|  ~ d \alpha \ ds
\tag{1.11}
\end{gathered}
\end{equation}
where $|J_1|=
\left|\begin{array}{cccc} \frac{\partial r}{\partial s} & \frac{\partial r}{\partial \alpha} \\
\frac{\partial \phi}{\partial s} & \frac{\partial \phi}{\partial \alpha} \end{array} \right|$
and similarly $|J_2|$ ~ (by Table 1 and  $\frac{\partial \theta}{\partial \lambda}=1$) are the corresponding Jacobians.
\end{theorem}
The complicated formulas above need numerical approximations by computer.
\subsection{Regular bounded periodic prism tilings and their space groups $\mathbf{pq2_1}$}
In \cite{Sz13-1} we have defined and described the regular prisms and prism tilings with a space group class $\Gamma=\mathbf{pq2_1}$ of $\SLR$.
These will be summarized in this section.
\begin{Definition}
Let $\cP^i$ be an infinite solid that is bounded by certain surfaces that can be
determined (in \cite{Sz13-1}) by ,,side fibre lines" passing through the
vertices of a regular $p$-gon $\cP^b$ lying in the base plane.
The images of solids $\cP^i$ by $\SLR$ isometries are called {\rm infinite regular $p$-sided prisms}.
Here regular means that the side surfaces are congruent to each other under rotations about a fiber
line (e.g. through the origin).
\end{Definition}
The common part of $\cP^i$ with the base plane is the {\it base figure} of $\cP^i$ that is denoted by $\cP$ and its vertices coincide
with the vertices of $\cP^b$, {\bf but $\cP$ is not assumed to be a polygon}.
\begin{Definition}
A {\rm bounded regular $p$-sided prism} is analogously defined if the face of the base figure $\cP$ and its translated copy $\cP^t$, under
a fibre translation by (1.2) and so (1.3), are also introduced.
The faces $\cP$ and $\cP^t$ are called {\rm cover faces}.
\end{Definition}
We consider regular prism tilings $\cT_p(q)$ by prisms $\cP_p(q)$ where $q$ pieces regularly meet
at each side edge by $q$-rotation.

The following theorem has been proved in \cite{Sz13-1}:
\begin{theorem}
There exist regular bounded not face-to-face prism tilings $\cT_p(q)$ in $\SLR$ for each $3 \le p \in \mathbb{N}$ where $\frac{2p}{p-2} < q \in \mathbb{N}$.
\end{theorem}
We assume that the prism $\cP_p(q)$ is a {\it topological polyhedron} having at each vertex
one $p$-gonal cover face (it is not a polygon at all) and two {\it skew quadrangles} which lie on certain side surfaces in the model.
Let $\cP_p(q)$ be one of the tiles of $\cT_p(q)$, $\cP^b$ is centered in the origin with vertices $A_1A_2A_3 \dots A_p$ in the base plane (Fig.~1 and 2).
It is clear that the side curves $c_{A_iA_{i+1}}$ $(i=1\dots p, ~ A_{p+1} \equiv A_1)$
of the base figure are derived from each other by $\frac{2\pi}{p}$ rotation about the vertical $x$ axis, so there are congruent in $\SLR$ sense.
The corresponding vertices $B_1B_2B_3 \dots B_p$ are generated by a fibre translation $\tau$ given by (1.3)
with parameter $0< \Phi\in \mathbb{R}$.
\begin{figure}[ht]
\centering
\includegraphics[width=12cm]{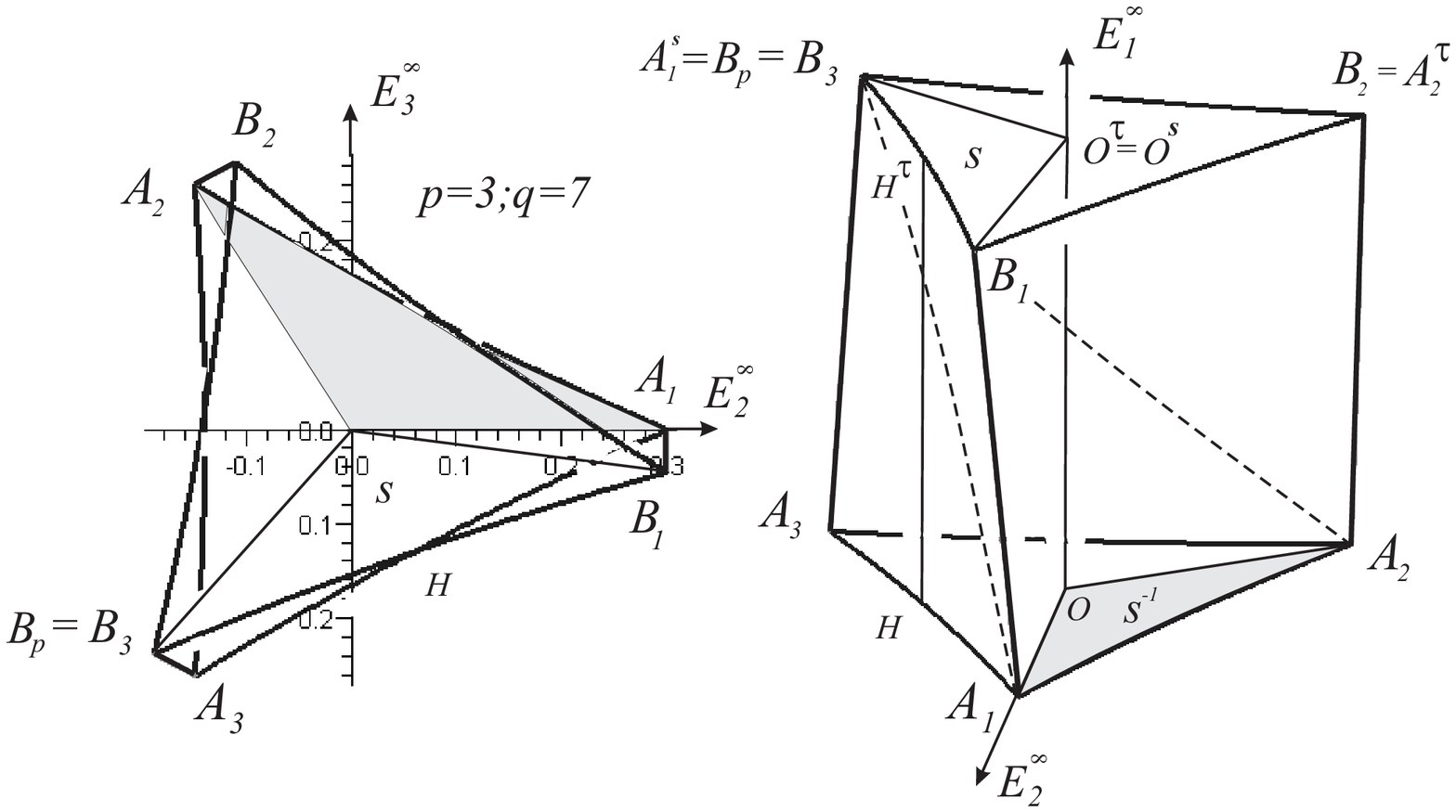}
\caption{The regular prism $\cP_p(q)$ and the fundamental domain of the space group ${\mathbf{pq2_1}}$}
\label{}
\end{figure}
The fibre lines through the vertices
$A_iB_i$ are denoted by $f_i, \ (i=1, \dots, p)$ and the fibre line through the "midpoint" $H$ of the curve $c_{A_1A_{p}}$ is denoted by
$f_0$. This $f_0$ will be a half-screw axis as follows below.

The tiling $\cT_p(q)$ is generated by a
discrete isometry group $\Gamma_p(q)=\mathbf{pq2_1}$ $\subset Isom(\SLR)$
which is given by its fundamental domain $A_1A_2O A_1^{\bs} A_2^{\bs} O^{\bs}$ a {\it topological polyhedron} and the group presentation
(see Fig.~1 and 4 for $p=3$ and \cite{Sz13-1} for details):
\begin{equation}
\begin{gathered}
\mathbf{pq2_1}=\{ \ba,\bb,\bs: \ba^p=\bb^q=\ba \bs \ba^{-1} \bs^{-1}= \bb \ba \bb \bs^{-1}=\mathbf{1} \}= \\
= \{ \ba,\bb: \ba^p=\bb^q=\ba \bb \ba \bb \ba^{-1} \bb^{-1} \ba^{-1} \bb^{-1}=\mathbf{1} \}. \tag{1.12}
\end{gathered}
\end{equation}
Here $\ba$ is a ${p}$-rotation about the fibre line through the origin ($x$ axis),  $\bb$ is a ${q}$-rotation about the fibre line trough
$A_1$ and $\bs=\bb \ba \bb$ is a screw motion~ $\bs:~ OA_1A_2 \rightarrow O^{\bs} B_p B_1$. All these can be obtained by formulas (1.5) and (1.6).
Then we get that $\ba\bb\ba\bb=\bb\ba\bb\ba=:\tau$ is a fibre translation. 
Then $\ba \bb$ is a $\mathbf{2_1}$ half-screw motion about
$f_0=HH^{\tau}$ (look at Fig.~1) that also determines the fibre tarnslation $\tau$ above. This group in (3.1) surprisingly occurred in \S ~ 6 of our paper \cite{MSzV} at double links
$K_{p,q}$.
The coordinates of the vertices $A_1A_2A_3 \dots A_p$ of the base figure and the corresponding vertices $B_1B_2B_3 \dots B_p$
of the cover face can be computed for all given parameters $p,q$ by
\begin{equation}
\tanh(OA_1)=b:=\sqrt{\frac{1-\tan{\frac{\pi}{p}} \tan{\frac{\pi}{q}}} {1+\tan{\frac{\pi}{q}} \tan{\frac{\pi}{q}}}}. \tag{1.13}
\end{equation}
\subsection{The volume of the bounded regular prisms}
The volume formula of a {\it sector-like} 3-dimensional domain $Vol(D(\Psi))$ can standardly be computed by the metric tensor $g_{ij}$ (see \cite{MSz14}).
in hyperboloid coordinates. This defined by the base figure $D$ lying in the base plane and by
fibre translation $\tau$ given by (1.3) with the height parameter $\Psi$.
\begin{theorem}
Suppose we are given a sector-like region $D$, so a continuous function $r = r(\theta)$
where the radius $r$ depends upon the polar angle $\theta$. The volume of domain $D(\Psi))$ is derived by the following integral:
\begin{equation}
\begin{gathered}
Vol(D(\Psi))=\int_{D}  \frac{1}{2}\sinh(2r(\theta)) {\mathrm{d}}r~ {\mathrm{d}} \theta  ~{\mathrm{d}} \psi = \\ = \int_0^{\Psi} \int_{\theta_1}^{\theta_2} \int_{0}^{r(\theta)}
\frac{1}{2}\sinh(2r(\theta))~ {\mathrm{d}}r~ {\mathrm{d}} \theta \ {\mathrm{d}} \psi =\Psi \int_{\theta_1}^{\theta_2}
\frac{1}{4}(\cosh(2r(\theta))-1)~ {\mathrm{d}} \theta.
\tag{1.14}
\end{gathered}
\end{equation}
\end{theorem}
$\cP_p(q)$ be an arbitrary bounded regular prism. We get the following
\begin{theorem}
The volume of the bounded regular prism $\cP_p(q)$ \Big($3 \le p \in \mathbb{N} $, $\frac{2p}{p-2} < q \in \mathbb{N}$\Big) can be computed by the following simple formula:
\begin{equation}
Vol(\cP_p(q))=Vol(D(p,q,\Psi)) \cdot p, \tag{1.15}
\end{equation}
where $Vol(D(p,q,\Psi))$ is the volume of the sector-like 3-dimensional domain that is given by the sector region $OA_1A_2 \subset \cP$ (see Fig.~1 and 3)
and by $\Psi$ the $\SLR$ height of the prism, depending on $p,q$.
\end{theorem}
\begin{figure}[ht]
\centering
\includegraphics[width=5cm]{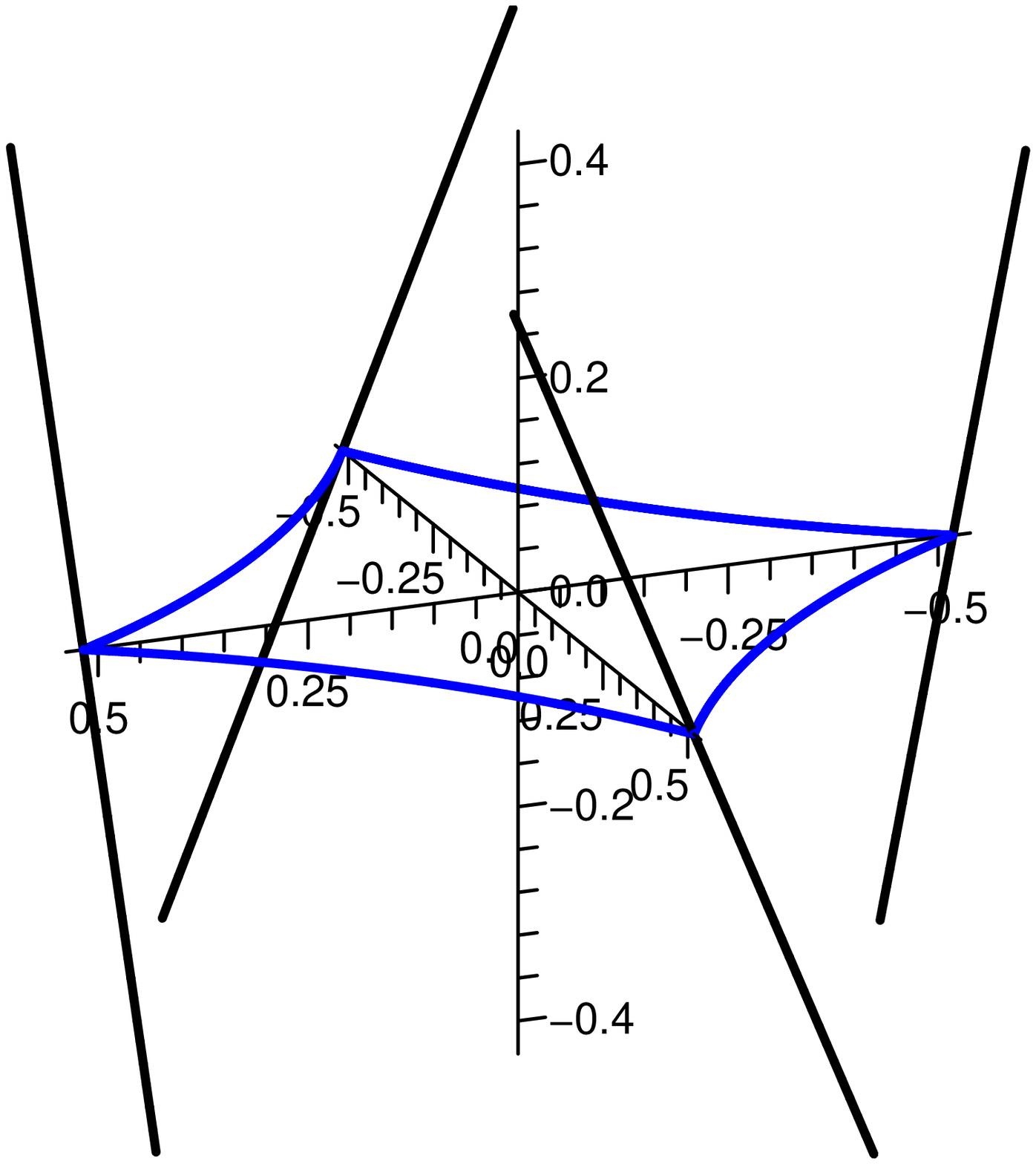} \includegraphics[width=6cm]{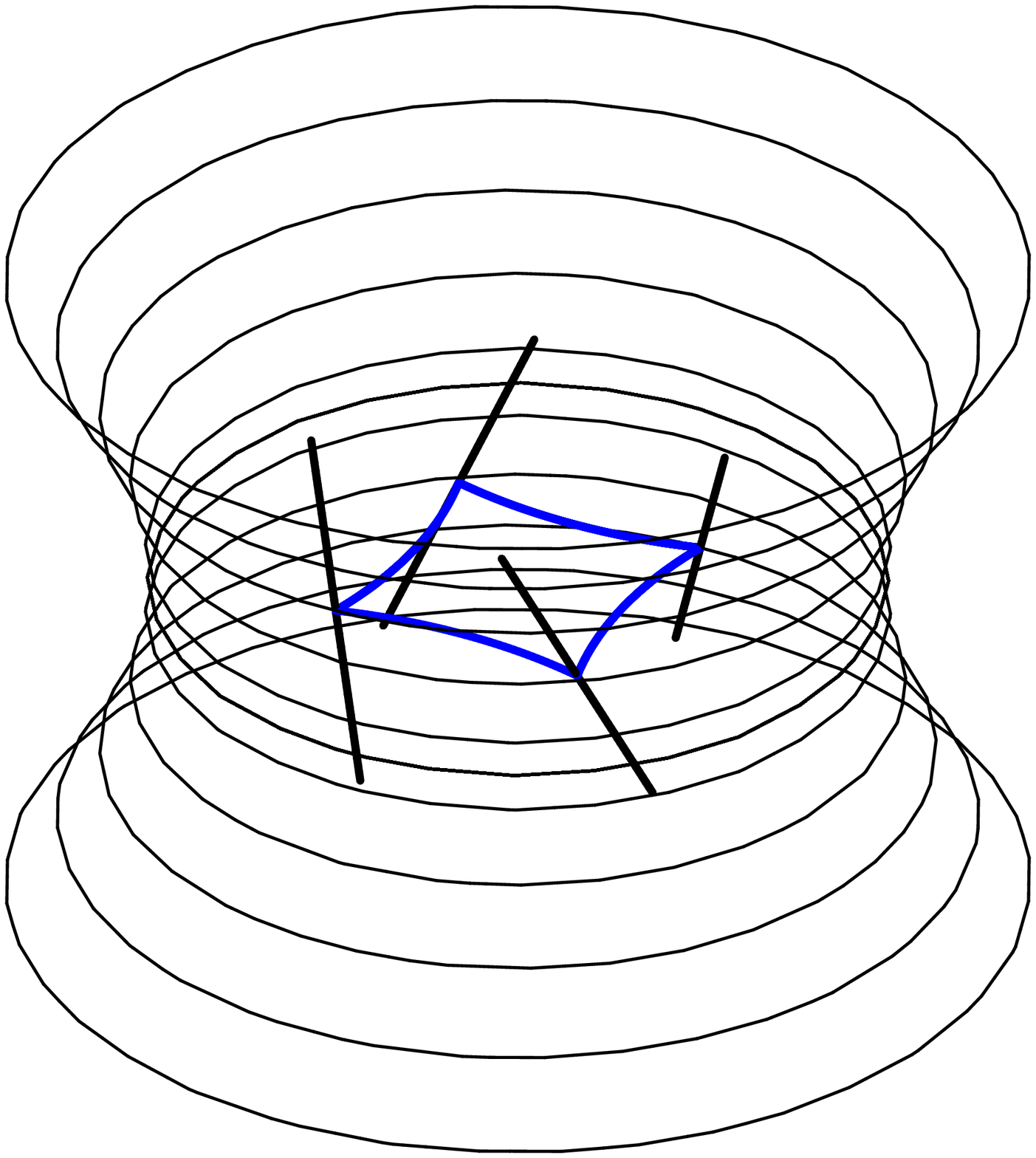}
\caption{Regular infinite 4-gonal prism $\cP_4^i(6)$ of the infinite regular prism tiling $\cT_4^i(6)$}
\label{}
\end{figure}
\section{Regular infinite prism tilings and non-periodic ball packings}
\subsection{Infinite regular prism tilings}
In this subsection we study the regular infinite prism tilings $\cT^i_p(q)$.
Let $\cT_p(q)$ be a regular prism tiling and let 
$\cP_p(q)$ be one of its tiles which is given by its base figure $\cP$ that is centered at the origin $K$ with vertices $G_1G_2G_3 \dots G_p$
in the base plane of the model and the corresponding vertices $A_1A_2$ $A_3 \dots A_p$ and $B_1B_2B_3 \dots B_p$
are generated by fibre translations $-\tau$ and $\tau$ given by (1.3)
with parameter $\Psi=\frac{\pi}{2}-\frac{\pi}{p}-\frac{\pi}{q}$. The images of the topological polyhedron $\cP_p(q)$ 
by the translations $\langle \tau \rangle $ form an infinite prism $\cP^i_p(q)$ (see Definitions 1.~4-5).
\begin{figure}[ht]
\centering
\includegraphics[width=12cm]{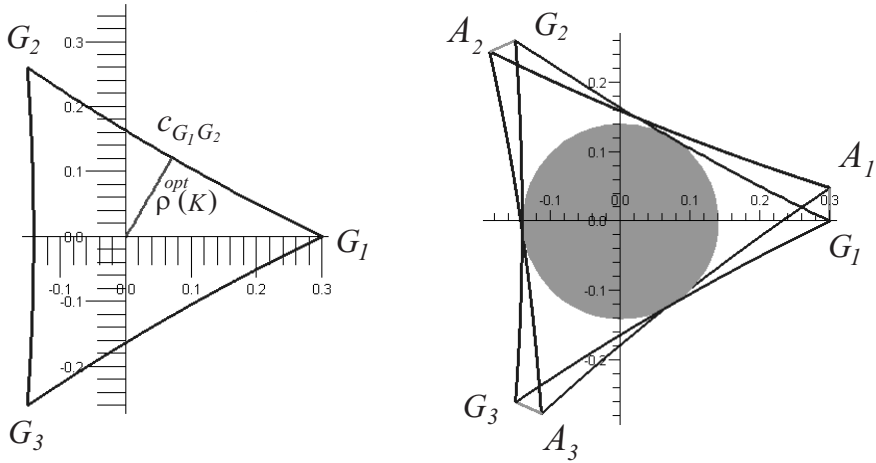} 
\caption{The maximal radius $\rho^{opt}(K)$ and the optimal half prism $A_1A_2A_3G_1G_2G_3$ with optimal half sphere for parameters $p=3$, $q=7$ with the
maximal radius}
\label{}
\end{figure}
By the constuction of the bounded prism tilings follows that rotations through $\omega=\frac{2\pi}{q}$ about the fibre lines
$f_i$ maps the corresponding side face onto the neighbouring one.
Therefore, we have got the following (see \cite{Sz13-1}):
\begin{theorem}
There exist regular infinite face-to-face prism tilings $\cT_p^i(q)$ for integer parameters $p,q$ where $3 \le p,~ \frac{2p}{p-2}<q$.
\end{theorem}
For example, we have described $\cP_4^i(6)$ with its base polygon in Fig.~2, where the parameter $b=\frac{\sqrt{6}-\sqrt{2}}{2}$.
\subsection{Non-periodic geodesic ball packings}
We consider a infinite regular prism tiling $\cT_p^i(q)$ and let $\cP_p^i(q)$ one of its tiles with base figure $\cP$ 
centered at the origin with vertices $G_1G_2 \dots G_p$
in the base plane of the model. 
Let $B_K^{opt}$ be the geodesic ball with center at the origin $K$ that touches the side surfaces of the infinite
regular prism $\cP_p^i(q)$. The radius of the ball $B_K^{opt}$ is denoted by $\rho^{opt}(K)$. Moreover,
we define the regular prism $\cP^{opt}_p(q)=A_1 A_2 \dots A_p B_1 B_2 \dots B_p$ with base figur $\cP$ and with cover faces
$A_1 A_2 \dots A_p$ and $B_1 B_2 \dots B_p$ touching $B_K^{opt}$. 
It is clear, that the height $h^{opt}_p(q)$ of $\cP^{opt}_p(q)$ is $2\rho^{opt}(K)$. 

The images of $\cP^{opt}_p(q)$ 
by the fibre traslations $\langle \tau \rangle$ where $h^{opt}_p(q)=|\tau|=2\rho^{opt}(K)$ cover the infinite regular prism $\cP_p^i(q)$ and 
by the structure of the infinte prism tilings follows that rotations through $\omega=\frac{2\pi}{q}$ about the fibre lines
$f_i$ maps the corresponding side face onto the neighbouring one and thus the images of $\cP^{opt}_p(q)$ fill the $\SLR$ space without overlap. These
tilings are denoted by $\cT^{n}_p(q)$.

{\it The height $h^{opt}_p(q)$ of the prism $\cP^{opt}_p(q)$ is not equal to $\pi-\frac{2\pi}{p}-\frac{2\pi}{q}$ so the corresponding 
regular prism tiling is non-periodic.}
We note here, that there are infinitely many non-periodic prism tilings derived from $\cT^{n}_p(q)$.

For the density of the packing it is sufficient to relate the volume of the optimal ball
to that of the solid $\cP^{opt}_p(q)$. The densitiy of the optimal ball packing of the prism tiling $\cT^{n}_p(q)$ ($3 \le p,~ \frac{2p}{p-2}<q$, 
integer parameters) can be computed by the following formula:  
\begin{equation}
\delta^{opt}_p(q):=\frac{Vol(B_K^{opt})}{Vol(\cP^{opt}_p(q))}. \notag
\end{equation}
In order to determine the optimal radius $\rho^{opt}(K)$ we will use the following Lemmas.
\begin{figure}[ht]
\centering
\includegraphics[width=8cm]{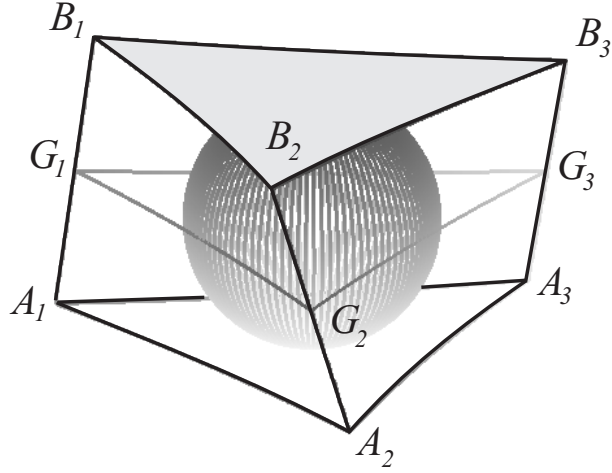} 
\caption{The optimal prism $A_1A_2A_3B_1B_2B_3$ with optimal sphere for parameters $p=3$, $q=7$ with the
maximal radius $\rho^{opt}(K)$}
\label{}
\end{figure}
The equation of the side curve $c_{G_1G_2}$ is derived as the foot points 
(see (1.3) and (1.8)) of the corresponding fibre lines ($3 \le p,~ \frac{2p}{p-2}<q$, where $p$ and $q$ are integer parameters):
\begin{lemma}
The parametric equation of the side curve $c_{G_1G_2}$ of the base figur $\cP$ is 
\begin{equation}
\scriptsize
\begin{gathered}
c_p^q(t)=\Bigg(0,~\sqrt{\sin \left( \frac{2\pi}{p}+\frac{2\pi}{q} \right)}\Bigg( t \cos \left( \frac{2\pi}{p} \right)\sin^2\left(\frac{\pi}{p}+\frac{\pi}{q}\right)-
\frac{t}{2}\sin\left(\frac{2\pi}{p}\right)\sin\left(\frac{2\pi}{p}+\frac{2\pi}{q}\right)+ \\
\sin^2\left(\frac{\pi}{p}+\frac{\pi}{q}\right)(1-t)+ t^2 \cos\left(\frac{\pi}{p}+\frac{\pi}{q}\right)
\cos\left(\frac{\pi}{p}-\frac{\pi}{q}\right)\Bigg)\Big/ \\ 
\Bigg({\sqrt{\left( \sin \left( {\frac {2\pi }{p}} \right) +
\sin \left( {\frac {2\pi }{q}} \right)  \right)}} \Big(\sin^2\left(\frac{\pi}{p}+\frac{\pi}{q}\right)+t^2 \cos^2\left(\frac{\pi}{p}+\frac{\pi}{q}\right)\Big) \Bigg),\\
t\sqrt{\sin \left( \frac{2\pi}{p}+\frac{2\pi}{q} \right)}\Bigg( \sin \left( \frac{2\pi}{p} \right)\sin^2\left(\frac{\pi}{p}+\frac{\pi}{q}\right)+
\frac{1}{2}\cos\left(\frac{2\pi}{p}\right)\sin\left(\frac{2\pi}{p}+\frac{2\pi}{q}\right)(1-t)+ \\
\cos\left(\frac{\pi}{p}+\frac{\pi}{q}\right)\Big(t \sin\left(\frac{2\pi}{p}\right)\cos\left(\frac{\pi}{p}+\frac{\pi}{q}\right)+
\sin\left(\frac{\pi}{p}+\frac{\pi}{q}\right)(t-1)\Big)\Bigg)\Big/ \\ 
\Bigg({\sqrt{\left( \sin \left( {\frac {2\pi }{p}} \right) +
\sin \left( {\frac {2\pi }{q}} \right)  \right)}} \Big(\sin^2\left(\frac{\pi}{p}+\frac{\pi}{q}\right)+t^2 
\cos^2\left(\frac{\pi}{p}+\frac{\pi}{q}\right)\Big) \Bigg), ~ ~ t \in [0,1].
\end{gathered} {\tag{\normalsize{2.1}}}
\end{equation}
\end{lemma}
The side curves $c_{G_iG_{i+1}}$ $(i=1\dots p, ~ G_{p+1} \equiv G_1)$
of the base figure are derived from each other by $\frac{2\pi}{p}$ rotation about the vertical $x$ axis, so there are congruent and 
their curvatures are equal in $\SLR$ sense.
Moreover, the above side curves are congruent also in Euclidean sense, therefore their curvatures are equal in Euclidean sense, as well. We obtain 
by the usual machinery of the differential geometry the next  
\begin{lemma}
The curvature $C_p(q)$ of the side curves $c_{G_iG_{i+1}}$ $(i=1\dots p, ~ G_{p+1} \equiv G_1)$ in the Euclidean sense is 
\begin{equation}
C_p(q)=\sqrt {\frac{\cos \left( {\frac {\pi }{p}}+{\frac {\pi }{q}} \right) 
\left( \sin \left( {\frac {2\pi }{p}} \right) +\sin \left( {
\frac {2\pi }{q}} \right)  \right)} { \sin \left( {\frac {\pi }{p}
}+{\frac {\pi }{q}}  \right)  \left( 1-\cos \left( {
\frac {2\pi }{p}} \right)  \right)}} \tag{2.2}
\end{equation}
therefore, the side curves $c_{G_iG_{i+1}}$ $(i=1\dots p, ~ G_{p+1} \equiv G_1)$ are Euclidean circular arcs of radius $r_p^q=\frac{1}{C_p(q)}$.
\end{lemma}
\begin{Remark}
\begin{enumerate}
\item It is easy to see, that the asymptotic behaviour of $C_p(q)$ is the following:
$
\lim_{q~\rightarrow \infty}(C_p(q))=\cot\left( \frac {\pi }{p} \right),~ ~ \lim_{p~\rightarrow \infty}(C_p(q))=\infty.
$
\item 
Given a point off of a line, if we drop a perpendicular to the above line from the given point, then $x$ is the distance 
along this perpendicular segment, and let $\phi=\Pi(x)$ is the least angle such that the line drawn 
through the point at that angle does not intersect the given line. 
The angle $\phi$ is the angle of parallelism. By the famous formel of J.~Bolyai follows, that $\log (\cot(\phi))=x$. Therefore, 
if we denote the distance of parallelism of the angle $\phi$ by $\Lambda(\phi)$ then 
$\log \Big( \lim_{q~\rightarrow \infty}(C_p(q))\Big)=\log \Big(\cot\left( \frac {\pi }{p} \right)\Big)=\Lambda\Big(\frac {\pi }{p}\Big)$.
\end{enumerate}
\end{Remark}
In the Table 1 we have collected some values of the radii of curvature $r_3^q$ of the side curve $c_{G_1G_{2}}$
of the base figur $\cP$. 
\medbreak
\centerline{\vbox{
\halign{\strut\vrule\quad \hfil $#$ \hfil\quad\vrule
&\quad \hfil $#$ \hfil\quad\vrule &\quad \hfil $#$ \hfil\quad\vrule &\quad \hfil $#$ \hfil\quad\vrule &\quad \hfil $#$ \hfil\quad\vrule
\cr
\noalign{\hrule}
\multispan5{\strut\vrule\hfill{Table 1} \hfill\vrule}%
\cr
\noalign{\hrule}
\noalign{\vskip2pt}
\noalign{\hrule}
(p, q) & (3,7) &  (3,8) & (3,10) & (3,1000) \cr
\noalign{\hrule}
C_p(q) & 0.286926 &  0.371579  & 0.453885  &  0.577339   \cr
\noalign{\hrule}
r_p^q & 3.485219 & 2.691215 & 2.203203 & 1.732085  \cr
\noalign{\hrule}
}}}
\medbreak
The maximal radius $\rho^{opt}(K)$ of the balls $B_K^{opt}$ can be determined using the above Lemmas  
for all possible parameters as the distance between the origin and $c_{G_1G_2}$. The volumes $Vol(B_K^{opt})$
can be computed by the Theorem 1.3 and the volumes of the prisms $\cP^{opt}_p(q)$ can be determined by the Theorem 1.4.

The above locally densest geodesic ball packings can be determined for all regular prism tilings $\cT^n_p(q)$ ($p,q$ as above). 
We have summarized in the following Tables some results to tilings $\cT_p^n(q)$. 
\medbreak
\centerline{\vbox{
\halign{\strut\vrule\quad \hfil $#$ \hfil\quad\vrule
&\quad \hfil $#$ \hfil\quad\vrule &\quad \hfil $#$ \hfil\quad\vrule &\quad \hfil $#$ \hfil\quad\vrule &\quad \hfil $#$ \hfil\quad\vrule
\cr
\noalign{\hrule}
\multispan5{\strut\vrule\hfill{Table 2} \hfill\vrule}%
\cr
\noalign{\hrule}
\noalign{\vskip2pt}
\noalign{\hrule}
(p, q)& \rho^{opt}(K) &  Vol(B_K^{opt}) & Vol(\cP^{opt}_p(q)) & \delta^{opt}_p(q) \cr
\noalign{\hrule}
\noalign{\hrule}
(3,7) &  0.141564 & 0.011963 & 0.031767 & 0.376592 \cr
\noalign{\hrule}
(3,8) &  0.181760 & 0.025431 & 0.071377 & 0.356287 \cr
\noalign{\hrule}
(3,10) & 0.219795 & 0.045198 & 0.138101 & 0.327281 \cr
\noalign{\hrule}
(3,1000) & 0.274648 & 0.088981 & 0.428828& 0.207499 \cr
\noalign{\hrule}
\vdots & \vdots & \vdots & \vdots & \vdots \cr
\noalign{\hrule}
(4,5) &  0.265319 & 0.080085 & 0.166705 & 0.480397 \cr
\noalign{\hrule}
(4,6) &  0.329239 & 0.154965 & 0.344779 & 0.449464 \cr
\noalign{\hrule}
(4,10) & 0.404230 & 0.292043 & 0.761956 & 0.383280 \cr
\noalign{\hrule}
(4,1000) & 0.440683 & 0.382228 & 1.378910 & 0.277196 \cr
\noalign{\hrule}
\vdots & \vdots & \vdots & \vdots & \vdots \cr
\noalign{\hrule}
(5,4) &  0.313435 & 0.133256 & 0.246171 & 0.541312 \cr
\noalign{\hrule}
(5,5) &  0.421241 & 0.332010 & 0.661684 & 0.501765 \cr
\noalign{\hrule}
(5,10) & 0.530638 & 0.686600 & 1.667047 & 0.411866 \cr
\noalign{\hrule}
(5,1000) & 0.562086 & 0.825191 & 2.639937 & 0.312580 \cr
\noalign{\hrule}
\vdots & \vdots & \vdots & \vdots & \vdots \cr
\noalign{\hrule}
(6,4) &  0.440687 & 0.382237 & 0.692229 & 0.552183 \cr
\noalign{\hrule}
(6,5) &  0.530638 & 0.686600 & 1.333638 & 0.514833 \cr
\noalign{\hrule}
(6,10) & 0.629251 & 1.188024 & 2.767592 & 0.429263 \cr
\noalign{\hrule}
(6,1000) & 0.658476 & 1.377893 & 4.124915 & 0.334042 \cr
\noalign{\hrule}
\vdots & \vdots & \vdots & \vdots & \vdots \cr
\noalign{\hrule}
(7,3) &  0.272637 & 0.087010 & 0.142753 & 0.609513 \cr
\noalign{\hrule}
(7,4) &  0.535202 & 0.705586 & 1.261041 & 0.559527 \cr
\noalign{\hrule}
(7,5) &  0.617496 & 1.117400 & 2.133913 & 0.523639 \cr
\noalign{\hrule}
(7,10) & 0.710652 & 1.772033 & 4.018646 & 0.440953 \cr
\noalign{\hrule}
(7,1000) & 0.738668 & 2.015812 & 5.785244 & 0.348440 \cr
\noalign{\hrule}
\vdots & \vdots & \vdots & \vdots & \vdots \cr
\noalign{\hrule}
(8,3) &  0.382143 & 0.245334 & 0.400179 & 0.613062 \cr
\noalign{\hrule}
(8,4) &  0.612113 & 1.086117 & 1.923010 & 0.564800 \cr
\noalign{\hrule}
(8,5) &  0.690221 & 1.608804 & 3.035751 & 0.529953 \cr
\noalign{\hrule}
(8,10) & 0.780165 & 2.422804 & 5.392115 & 0.449324 \cr
\noalign{\hrule}
(8,1000) & 0.807443 & 2.722797 & 7.589676 & 0.358750 \cr
\noalign{\hrule}
\vdots & \vdots & \vdots & \vdots & \vdots \cr
\noalign{\hrule}
}}}
\medbreak
\centerline{\vbox{
\halign{\strut\vrule\quad \hfil $#$ \hfil\quad\vrule
&\quad \hfil $#$ \hfil\quad\vrule &\quad \hfil $#$ \hfil\quad\vrule &\quad \hfil $#$ \hfil\quad\vrule &\quad \hfil $#$ \hfil\quad\vrule
\cr
\noalign{\hrule}
\multispan5{\strut\vrule\hfill{Table 3} \hfill\vrule}%
\cr
\noalign{\hrule}
\noalign{\vskip2pt}
\noalign{\hrule}
(p, q)& \rho^{opt}(K) &  Vol(B_K^{opt}) & Vol(\cP^{opt}_p(q)) & \delta^{opt}_p(q) \cr
\noalign{\hrule}
\noalign{\hrule}
(10,3) &  0.530638 & 0.686600 & 1.111365 & 0.617799 \cr
\noalign{\hrule}
\vdots & \vdots & \vdots & \vdots & \vdots \cr
\noalign{\hrule}
(20,3) &  0.914848 & 4.195479 & 6.706186 & 0.625613 \cr
\noalign{\hrule}
(20,4) &  1.094612 & 8.023914 & 13.755306 & 0.583332 \cr
\noalign{\hrule}
(20,5) &  1.163424 & 10.092704 & 18.275027 & 0.552268 \cr
\noalign{\hrule}
(20,10) & 1.245625 & 13.132701 & 27.392724 & 0.479423 \cr
\noalign{\hrule}
(20,1000) & 1.271043 & 14.216772 & 35.858024 & 0.396474 \cr
\noalign{\hrule}
\vdots & \vdots & \vdots & \vdots & \vdots \cr
\noalign{\hrule}
(28,3) & 1.088398 & 7.855861 & 12.537440 & 0.626592 \cr
\noalign{\hrule}
\mathbf{(29,3)} & \mathbf{1.106311} & \mathbf{8.348310} & \mathbf{13.323054} & \mathbf{0.626606} \cr
\noalign{\hrule}
(30,3) & 1.123593 & 8.847342 & 14.119487 & 0.626605 \cr
\noalign{\hrule}
\vdots & \vdots & \vdots & \vdots & \vdots \cr
\noalign{\hrule}
(35,3) & 1.201914 & 11.432334 & 18.250297 & 0.626419 \cr
\noalign{\hrule}
\vdots & \vdots & \vdots & \vdots & \vdots \cr
\noalign{\hrule}
(40,3) &  1.269482 & 14.148085 & 22.599777 & 0.626028 \cr
\noalign{\hrule}
\vdots & \vdots & \vdots & \vdots & \vdots \cr
\noalign{\hrule}
(52,3) &  1.401728 & 21.089811 & 33.761388 & 0.624673 \cr
\noalign{\hrule}
\vdots & \vdots & \vdots & \vdots & \vdots \cr
\noalign{\hrule}
(72,3) &  1.565173 & 33.642710 & 54.088487 & 0.621994 \cr
\noalign{\hrule}
}}}
\medbreak
\begin{Remark}
\begin{enumerate}
\item The best density that we found $\approx 0.626606$ for parameters $p=29, q=3$ that is larger  
then the maximal density of the corresponding periodical geodesic ball packings under the groups $\mathbf{pq2_1}$. 
\item The problems of finding the densest geodesic and translation ball packings in the Thurston gemetries are timely 
(see e.g. \cite{MSz12}, \cite{Sz07}, \cite{Sz10-2}, \cite{Sz10-3}, \cite{Sz12-1}).
\end{enumerate}
\end{Remark}
%
%{\bf{Acknowledgement:}}

 \end{document}